\documentclass{article}
\usepackage{amssymb,amsmath}

\newtheorem{theorem}{Theorem}
\newtheorem{lemma}{Lemma}
\newtheorem{problem}{Problem}
\newcommand{\bt}{\begin{theorem}}
\newcommand{\et}{\end{theorem}}
\newcommand{\bl}{\begin{lemma}}
\newcommand{\el}{\end{lemma}}
\newcommand{\bp}{\begin{problem}}
\newcommand{\ep}{\end{problem}}
\newcommand{\pf}{{\bf Proof}.\ }
\newcommand{\eop}{ $\square$ \vspace{.8cm}}
\newcommand{\beal}{\begin{align*}}
\newcommand{\enal}{\end{align*}}
\newcommand{\beq}{\begin{equation}}
\newcommand{\eeq}{\end{equation}}
\newcommand{\benum}{\begin{enumerate}}
\newcommand{\eenum}{\end{enumerate}}
\newcommand{\ba}{\begin{array}}
\newcommand{\ea}{\end{array}}

\newcommand{\N}{\ensuremath{ \mathbf N }}
\newcommand{\FF}{\ensuremath{\mathcal F}}
\newcommand{\GG}{\ensuremath{\mathcal G}}

\newcommand{\supp}{\text{supp}}

\newcommand{\pol}{$\mathcal{F} = \{f_n(q)\}_{n=1}^{\infty}$}

\newcommand{\polg}{$\mathcal{G} = \{g_n(q)\}_{n=1}^{\infty}$}

\newcommand{\tfe}{ the functional equation~(\ref{qps:fe})}

\begin{document}
\title{Formal power series arising from \\
multiplication of quantum integers}
\author{Melvyn B. Nathanson\thanks{This work was supported
in part by grants from the NSA Mathematical Sciences Program
and the PSC-CUNY Research Award Program.}\\
Department of Mathematics\\
Lehman College (CUNY)\\
Bronx, New York 10468\\
Email: nathansn@alpha.lehman.cuny.edu}
\maketitle

\begin{abstract}
For the quantum integer $[n]_q = 1+q+q^2+\cdots + q^{n-1}$
there is a natural polynomial multiplication
such that $[mn]_q = [m]_q \otimes_q [n]_q$.  This multiplication is described by 
the functional equation $f_{mn}(q) = f_m(q) f_n(q^m)$, defined on a given sequence 
\pol\ of polynomials such that $f_n(0)=1$ for all $n$.  
If \pol\ is a solution of the functional equation, 
then there exists a formal power series $F(q)$ such that the sequence 
$\{f_n(q)\}_{n=1}^{\infty}$ converges to $F(q).$

Quantum mulitplication suggests the functional equation
\[
f(q)F(q^m) = F(q),
\]
where $f(q)$ is a fixed polynomial or formal power series 
with constant term $f(0)=1$,
and $F(q)=1+\sum_{k=1}^{\infty}b_kq^k$ is a formal power series.
It is proved that this functional equation has a unique solution 
$F(q)$ for every polynomial or formal power series $f(q)$.
If the degree of $f(q)$ is at most $m-1$,
then there is an explicit formula for the coefficients $b_k$ of $F(q)$ 
in terms of the coefficients of $f(q)$ and the $m$-adic representation of $k$.

The paper also contains a review of convergence properties of formal power series 
with coefficients in an arbitrary field or integeral domain.
\end{abstract}

\section{Quantum multiplication}
Nathanson~\cite{nath03b,nath03d} introduced the functional equation 
for multiplication of quantum integers as follows.

Let \N$ = \{1,2,3,\ldots\}$ denote the set of positive integers
and $\N_0 = \N \cup \{0\}$ the set of nonnegative integers.
For $n \in \N$, the polynomial
\[
[n]_q = 1 + q + q^2 + \cdots + q^{n-1}
\]
is called the {\em quantum integer} $n.$
With the usual multiplication of polynomials, however, we observe that
$[m]_q [n]_q \neq [mn]_q$ for all $m\geq 2$ and $n\geq 2$.
We would like to define a polynomial multiplication such that
the product of the quantum integers $[m]_q$ and $[n]_q$
is $[mn]_q$.

Let $\mathcal{F} = \{f_n(q)\}_{n=1}^{\infty}$ be a sequence of polynomials
with coefficients in a field.
We define a multiplication operation on the polynomials in $\mathcal{F}$ by
\[
f_m(q)\otimes_q f_n(q) = f_m(q)f_n(q^m).
\]
If $f_n(q) = [n]_q$ is the $n$-th quantum integer, then
\begin{align*}
[m]_q\otimes_q [n]_q
& = [m]_q[n]_{q^m} \\
& = \left(1+q+q^2+\cdots + q^{m-1}\right)\left(1+q^m+q^{2m} + \cdots + q^{m(n-1)}\right) \\
& = 1+q+q^2+\cdots + q^{mn-1} \\
& = [mn]_q
\end{align*}
for all positive integers $m$ and $n$,
and so the sequence $\mathcal{F} = \{ [n]_q\}_{n=1}^{\infty}$ 
satisfies the functional equation
\begin{equation}       \label{qps:fe}
f_{mn}(q) = f_m(q)\otimes_q f_n(q).
\eeq
This is called the {\em functional equation for quantum multiplication}.

Many other sequences of polynomials also satisfy this functional equation.
For example, for every positive integer $t$, the sequence of polynomials
$\{q^{t(n-1)}\}_{n=1}^{\infty}$ satisfies~(\ref{qps:fe}).

\bp    \label{qps:problem1}
Determine all sequences of polynomials \pol\
that satisfy the functional equation
\[
f_{mn}(q) = f_m(q)\otimes_q f_n(q).
\]
\ep

The sequence of \pol\ is called {\em nonzero} if $f_n(q) \neq 0$ for some integer $n$.
The functional equation implies that
\[
f_1(q) = f_1(q)^2,
\]
and so $f_1(q) = 0$ or 1.  
If $f_1(q) = 0,$ then
\[
f_n(q) = f_{1\cdot n}(q) = f_1(q)f_n(q) = 0
\]
for all $n \in \N.$
Therefore, a solution \pol\ of~({\ref{qps:fe}) is nonzero 
if and only if $f_1(q) = 1.$

The {\em support} of $\mathcal{F}$ is the set
\[
\supp(\mathcal{F}) = \{n \in \N: f_n(q) \neq 0\}.
\]
If a nonzero sequence of polynomials \pol\ is a solution of the functional
equation~(\ref{qps:fe}), then $\supp(\mathcal{F})$ is a multiplicative 
subsemigroup of the positive integers.
Nathanson~\cite{nath03b} proved that $\supp(\mathcal{F})$ is a semigroup of the form $S(P)$,
where $P$ is a set of prime numbers and $S(P)$ is the semigroup of positive integers
generated by $P$.

Every nonzero polynomial $f(q)$ with coefficients in a field can be 
written uniquely in the form 
\beq              \label{qps:polyform}
f(q) = aq^{v(f)}g(q),
\eeq
where $a\neq 0$, $v(f) \in \N_0$, and $g(q)$ is a polynomial with constant term $g(0) = 1$.
Let $\mathcal{F} = \{f_n(q)\}_{n=1}^{\infty}$ be a nonzero sequence of polynomials
that is a solution of the functional equation~(\ref{qps:fe}).
We represent each nonzero polynomial $f_n(q)$ in the form
\[
f_n(q) = a(n)q^{v(f_n)}g_n(q),
\]
where $a(n) \neq 0,$ $v(f_n)\in \N_0,$ and $g_n(0) = 1$.
We define $g_n(q) = 0$ and $a(n) = 0$ for all $n \in \N_0\setminus \supp(\mathcal{F}).$

For all $m,n \in \supp(\mathcal{F}),$
we have
\begin{align*}
a(mn)q^{v(f_{mn})}g_{mn}(q)
& = f_{mn}(q) \\
& = f_{m}(q) f_{n}(q^m) \\
& = a(m)q^{v(f_{m})}g_{m}(q)a(n)q^{mv(f_{n})}g_{n}(q^m)  \\
& = a(m)a(n) q^{v(f_m)+mv(f_n)}g_{m}(q)g_{n}(q^m).  
\end{align*}
Since
\[
g_{mn}(0) = g_{m}(0)g_{n}(0) = 1,
\]
it follows that 
\[
a(mn) = a(m)a(n)
\]
for all $m,n \in \supp(\FF)$, and so $a$ is a completely multiplicative 
arithmetic function on the semigroup $\supp(\FF)$
with $a(n) \neq 0$ for all $n \in \supp(\FF).$
Similarly, 
\[
q^{v(f_{mn})} = q^{v(f_m)+mv(f_n)},
\]
and so
\[
v(f_m)+mv(f_n) = v(f_{mn}) = v(f_{nm}) = v(f_n)+nv(f_m)
\]
for all integers $m$ and $n$ in $\supp(\FF)$.
This implies that there exists a nonnegative rational number $t$
such that 
\[
\frac{v(f_m)}{m-1} = \frac{v(f_n)}{n-1} = t
\]
for all $m,n \in \supp(\FF)\setminus\{1\}$.  Moreover,
\[
t(n-1) \in \N_0 
\]
and 
\[
q^{v(f_n)} = q^{t(n-1)}
\]
for all $n \in \supp(\FF)$.

Finally, we see that $\GG = \{g_n(q)\}_{n=1}^{\infty}$
is also a solution of the functional equation for quantum multiplication
with $\supp(\GG) = \supp(\FF)$, 
and $g_n(0) = 1$ for all $n \in \supp(\GG)$.
It follows that to solve Problem~\ref{qps:problem1} 
it suffices to classify solutions \pol\ of~(\ref{qps:fe})
such that $f_n(0) = 1$ for all $n \in \supp(\FF).$

Let $\mathcal{F} = \{f_n(q)\}_{n=1}^{\infty}$
be a solution of~(\ref{qps:fe}).
The operation $\otimes_q$ is commutative on $\mathcal{F}$ since
\[
f_m(q)\otimes_q f_n(q) = f_{mn}(q) = f_{nm}(q) = f_n(q)\otimes_q f_m(q).
\]
Equivalently,
\beq                 \label{qps:mnfe}
f_m(q)f_n(q^m) = f_n(q)f_m(q^n)
\eeq
for all positive integers $m$ and $n$.
Moreover, 
\[
\deg(f_m)+m\deg(f_n) = \deg(f_n) + n\deg(f_m),
\]
and so $\deg(f_m) = t(m-1)$ and $\deg(f_n) = t(n-1)$ for some nonnegative
rational number $t$.
This suggests the following problem.

\bp            \label{qps:problem2}
Let $m$ and $n$ be positive integers.
For a fixed polynomial $f_m(q)$, determine all polynomials $f_n(q)$ 
that satisfy the functional equation
\[
f_m(q)f_n(q^m) = f_n(q)f_m(q^n).
\]
\ep

Let $f_m(q) = \sum_{i=0}^{t(m-1)}a_iq^i$ and
$f_n(x) = \sum_{j=0}^{t(n-1)}b_jq^j$ be polynomials 
of degrees $t(m-1)$ and $t(n-1)$, respectively, with constant terms
$a_0 = b_0 = 1$.
If $f_m(q)$ and $f_n(q)$ satisfy equation~(\ref{qps:mnfe}),
then
\[
f_m(q) \sum_{j=0}^{t(n-1)}b_j q^{mj} 
= \left( \sum_{j=0}^{t(n-1)}b_j q^{j} \right)
\left(1 + a_1q^n +\cdots + a_{t(m-1)}q^{nt(m-1)}\right).
\]
Letting $n$ tend to infinity, we obtain the functional equation
\[
f_m(q) \sum_{j=0}^{\infty}b_j q^{mj} = \sum_{j=0}^{\infty}b_j q^{j},
\]
or, equivalently,
\[
f_m(q)F(q^m) = F(q),
\]
where 
\[
F(q) = \sum_{k=0}^{\infty}b_k q^k
\]
is a formal power series with constant term $b_0 = 1.$

\bp               \label{qps:problem3}
Let $m$ be positive integer and $f(q)$ a polynomial or formal power series 
with constant term $f(0)=1.$
Determine all formal power series $F(q)$ that satisfy the functional equation
\beq             \label{qps:psfe}
f(q)F(q^m) = F(q).
\eeq
\ep

Problems~\ref{qps:problem1} and~\ref{qps:problem2} are unsolved.
Problem~\ref{qps:problem3}, however, will be solved in Section~\ref{qps:section4},
where we prove that for every polynomial $f(q)$, and, more generally, 
for every formal power series $f(q)$
with constant term $f(0) = 1$, there is a unique formal power series 
$F(q)$ that is a solution of the functional equation
\[
f(q)F(q^m) = F(q).
\]
Moreover, we shall explicitly construct the coefficients of the formal power series 
$F(q)$ when $f(q)$ is a polynomial of degree at most $m-1.$

\section{Convergence of formal power series}\label{qps:sec:fps}
We review here some elementary convergence properties of formal power series
with coefficients in a ring.

\subsection{The ring of formal powers series}

A {\em formal power series} in the variable $q$
with coefficients in a ring $R$ is an expression of the form
\[
\sum_{n=0}^{\infty} a_n q^n = a_0 + a_1q = a_2q^2 + \cdots,
\]
where the coefficients $a_n$ are elements of $R$.
The coefficient $a_0$ is called the {\em constant term} of the power series.
We denote by $R[[q]]$ the set of all formal power series
with coefficients in $R$.
The sum and product of formal power series are defined, as usual, by
\[
\sum_{n=0}^{\infty} a_n q^n + \sum_{n=0}^{\infty} b_n q^n 
= \sum_{n=0}^{\infty} (a_n + b_n) q^n
\]
and
\[
\sum_{n=0}^{\infty} a_n q^n \cdot \sum_{n=0}^{\infty} b_n q^n  
= \sum_{n=0}^{\infty} c_n q^n,
\]
where
\[
c_n = \sum_{i=0}^{n}a_ib_{n-i}.
\]
With these operations of addition and multiplication,
$R[[q]]$ is a ring.

In this paper we shall always assume that $R$ is an integral domain.  
Then the ring of formal power series $R[[q]]$ is also an integral domain.

A formal power series $f(q)$ is {\em invertible}
if there exists a formal power series $g(q)$ such that $f(q) g(q) = 1$.
If $f(q) \in R[[q]]$ is invertible, then there is a unique 
formal power series $g(q)$ such that  $f(q) g(q) = 1$.
The inverse $g(q)$ is denoted $f(q)^{-1}$.
For example, $(1-q)^{-1} = \sum_{n=0}^{\infty} q^n$.

\bt  \label{qps:theorem:invert}
The formal power series $f(q) = \sum_{n=0}^{\infty} a_n q^n $
is invertible in $R[[q]]$ if and only if the coefficient $a_0$
is invertible in the ring $R$.
\et

\pf
If the power series $\sum_{n=0}^{\infty} b_nq^n$ is the 
multiplicative inverse of $f(q)$, then $a_0b_0 = 1$
and so $a_0$ is invertible in $R$.

Conversely, suppose that $a_0$ is invertible in $R$.
The formal power series $g(q) = \sum_{n=0}^{\infty}b_nq^n$
is a solution of the equation $f(q)g(q) = 1$ 
if and only if the coefficients $b_n$ satisfy the identities
\[
a_0b_0 = 1
\]
and
\[
\sum_{i=0}^n a_i b_{n-i}
= a_0b_n + a_1b_{n-1} + \cdots + a_nb_0 = 0 \qquad\text{for all $n \geq 1$.}
\]
These conditions are equivalent to
\[
b_0 = a_0^{-1}
\]
and
\[
b_n = -a_0^{-1} \sum_{i=1}^n a_i b_{n-i}\qquad\text{for all $n \geq 1$.}
\]
We can now construct inductively the unique sequence $\{b_n\}_{n=1}^{\infty}$
of elements of the ring $R$ such that $g(q) = \sum_{n=0}^{\infty} b_nq^n = f(q)^{-1}$.
This completes the proof.
\eop

Let $f(q) = \sum_{n=0}^{\infty} a_nq^n$ be a formal power series.
We define the {\em valuation} 
\[
v: R[[q]] \rightarrow \N_0 \cup \{\infty\}
\]
as follows:
If $f(q) \neq 0$, then $v(f)$ is the smallest integer $n$ such that
$a_n \neq 0$.  If $f(q) = 0$, then we set $v(f) = \infty$.
 
Let $f(q) = \sum_{n=0}^{\infty} a_n q^n$ and
$g(q) = \sum_{n=0}^{\infty} b_n q^n$ be formal power series.
We write
\[
f(q) \equiv g(q) \pmod{q^N}
\]
if $a_n = b_n$ for $n = 0,1,\ldots, N-1$.
This is equivalent to the inequality $v(f-g) \geq N$.

\bt               \label{qps:theorem:val} 
Let $R$ be an integral domain, and let $f,g \in R[[q]]$.  The valuation 
$v: R[[q]] \rightarrow \N_0 \cup \{\infty\}$ 
satisfies the following properties:
\benum
\item[(i)]
\[
v(-f) = v(f)
\]
\item[(ii)]
\[
v(fg) = v(f) + v(g)
\]
\item[(iii)]
\[
v(f \pm g) \geq \min(v(f),v(g))
\]
\item[(iv)]
\[
v(f \pm g) =  \min(v(f),v(g)) \qquad\text{if $v(f) \neq v(g)$.}
\]
\item[(v)]
\[
v(f(q^k)-f(0)) \geq k.
\]
\eenum
\et

\pf
Let $v(f) = s$ and $v(g) = t$.
Then $f(q) = \sum_{n=s}^{\infty} a_nq^n$ and
$g(q) = \sum_{n=t}^{\infty} b_nq^n$,
where $a_s$ and $b_{t}$ are nonzero elements of the ring $R$.
Then $-a_s \neq 0$, and so $v(-f) = s$.  This proves~(i).

Statement~(ii) is a consequence of the identity
\[
f(q)g(q) = a_s b_{t} q^{s + t} + \sum_{n=s+t+1}^{\infty} c_nq^n,
\]
where $c_n \in R$ for $n \geq s+t + 1$.
Since $R$ is an integral domain, it follows that $a_sb_{t} \neq 0$,
and so $v(fg) = s+t = v(f)+v(g)$.

To prove~(iii) and~(iv), we observe that if $s = t$, then
\[
f(q) \pm g(q) = \sum_{n=s}^{\infty} (a_n \pm b_n)q^n
\]
and $v(f \pm g) \geq s = \min(v(f),v(g))$.
If $s < t$, then
\[
f(q) \pm g(q) = \sum_{n=s}^{t -1} a_nq^n + \sum_{n=t}^{\infty} (a_n \pm b_n)q^n
\]
and so $v(f \pm g) = s = \min(v(f),v(g))$.

Finally, if $s=0,$ then $f(q)-f(0) =  \sum_{n=r}^{\infty} a_nq^n$,
where $v(f-f(0)) = r \geq 1$ and $a_r \neq 0.$
Then 
\[
f(q^k)-f(0) = a_rq^{kr} + \sum_{n=r+1}^{\infty} a_nq^{kn},
\]
and
\[
v(f(q^k)-f(0) = kr \geq k.
\]
If $s \geq 1,$ then $f(0) = 0$ and
\[
f(q^k)-f(0) = f(q^k)=  a_sq^{ks} + \sum_{n=s+1}^{\infty} a_nq^{kn},
\]
and
\[
v(f(q^k)-f(0) = ks \geq k.
\]
This proves~(v).
\eop

We introduce a topology on the ring of formal power series as follows.
Let $\{f_k(q)\}_{k=1}^{\infty}$ be a sequence
of formal power series in $R[[q]]$, where
\[
f_k(q) = \sum_{n=0}^{\infty} a_{k,n}q^n.
\]
Let
\[
f(q) = \sum_{n=0}^{\infty} a_{n}q^n \in R[[q]].
\]
The sequence $\{f_k(q)\}_{k=1}^{\infty}$ {\em converges} to $f(q)$,
that is,
\[
\lim_{k\rightarrow\infty} f_k(q) = f(q),
\]
if for every $n \geq 0$ there exists an integer $k_0(n)$
such that $a_{k,n} = a_n$ for all $k \geq k_0(n)$.
For example, the sequence $\{[k]_q\}_{k=1}^{\infty}$
of quantum integers converges to $(1-q)^{-1}$.

Similarly, if 
\[
f_{\mathbf{k}}(q) = \sum_{n=0}^{\infty} a_{\mathbf{k},n}q^n.
\]
is a formal power series for all 
\[
\mathbf{k} = (k_1,\ldots,k_d) \in \N_0^d,
\]
where $\N_0^d$ is the set of all $d$-tuples of nonnegative integers,
then
\[
\lim_{{\mathbf{k}}\rightarrow\infty} f_{\mathbf{k}}(q) 
\lim_{k_1,\ldots,k_d\rightarrow\infty} f_{k_1,\ldots,k_d}(q) = f(q),
\]
if for every $n \geq 0$ there exists an integer $k_0(n)$
such that $a_{{\mathbf{k}},n} = a_n$ for all ${\mathbf{k}} \in \N_0^d$
with $k_i \geq k_0(n)$ for $i = 1,\ldots,d.$

We consider convergence in the case $d = 2$ in Theorem~\ref{qps:theorem:convergence1}.

\bt              \label{qps:theorem:vconverge}
Let $\{f_k(q)\}_{k=1}^{\infty}$ be a sequence
of formal power series in $R[[q]]$.
The sequence $\{f_k(q)\}_{k=1}^{\infty}$ converges if and only if
\[
\lim_{k,\ell\rightarrow\infty} v(f_k-f_{\ell}) = \infty.
\]
Moreover, if $f(q) \in R[[q]]$, then
\[
\lim_{k\rightarrow\infty} f_k(q) = f(q)
\]
if and only if
\[
\lim_{k\rightarrow\infty} v(f-f_k) = \infty.
\]
\et

\pf
This follows immediately from the definitions of the valuation $v$
and convergence in the ring of formal power series.
\eop

\bt
If
\[
\lim_{k\rightarrow\infty} f_k(q) = f(q)
\qquad\text{and}\qquad
\lim_{k\rightarrow\infty} g_k(q) = g(q),
\]
then
\beq             \label{qps:add}
\lim_{k\rightarrow\infty} \left( f_k(q) + g_k(q) \right) = f(q) + g(q)
\eeq
and
\beq             \label{qps:mult}
\lim_{k\rightarrow\infty} f_k(q) g_k(q) = f(q)g(q).
\eeq
If $f_k(q)$ is invertible for all $k \geq 1$,
then $f(q)$ is invertible
and
\beq             \label{qps:inv}
\lim_{k\rightarrow\infty} f_k(q)^{-1} = f(q)^{-1}.
\eeq
\et

\pf
Statements~(\ref{qps:add}) and~(\ref{qps:mult}) are straightforward verifications.

To prove~(\ref{qps:inv}), let $a_{k,0}$ and $a_0$ be the constant terms 
of the formal power series $f_k(q)$ and $f(q)$, respectively, for all $k \geq 1.$ 
Since $f_k(q)$ is invertible, it follows that $a_{k,0}$ is a unit in $R$.
Since $a_0 = a_{k,0}$ for all sufficiently large $k$, it follows that
$a_0$ is a unit, and so $f(q)$ is invertible by Theorem~\ref{qps:theorem:invert}.
Applying Theorem~\ref{qps:theorem:val} to the identity
\[
f(q)^{-1} - f_k(q)^{-1} = (f_k(q) - f(q)) f(q)^{-1} f_k(q)^{-1},
\]
we obtain
\[
v(f^{-1} - f_k^{-1}) = v(f_k - f)+ v(f^{-1}) + v(f_k^{-1})
\geq v(f_k - f).
\]
Since $f_k(q)$ converges to $f(q)$, it follows that $v(f_k - f)$
tends to infinity, hence $v(f^{-1}-f_k^{-1})$ also tends to infinity
and $f_k(q)^{-1}$ converges to $f(q)^{-1}$.
This completes the proof.
\eop

\subsection{Infinite series and infinite products}

We define convergence of the {\em infinite series} \index{infinite series}
\[
\sum_{k=0}^{\infty} f_k(q)
\]
of formal power series as follows.
The {\em $m$--th partial sum}\index{partial sum}
\[
s_m(q) = \sum_{k=0}^{m-1} f_k(q)
\]
is simply a finite sum of formal power series.
If the sequence $\{s_m(q)\}_{m=0}^{\infty}$ converges to $f(q) \in R[[q]]$,
then we write that the infinite series converges to $f(q)$, that is,
\[
\sum_{k=0}^{\infty} f_k(q) = \lim_{m\rightarrow\infty} s_m(q) = f(q).
\]
For example, if $f_k(q) = \sum_{n=k}^{\infty}q^n$ for $k = 0,1,2,\ldots$,
then
\[
\sum_{k=0}^{\infty}f_k(q) = \sum_{n=0}^{\infty} (n+1)q^n.
\]

\bt               \label{qps:theorem:infiniteseries}
Let $\{f_k(q)\}_{k=0}^{\infty}$ be a sequence of formal power series.
The infinite series $\sum_{k=0}^{\infty}f_k(q)$ converges
if and only if $\lim_{k\rightarrow\infty} v(f_k) = \infty$.
\et

\pf
Let $f_k(q) = \sum_{n=0}^{\infty}a_{k,n}q^n$.
If $\lim_{k\rightarrow\infty} v(f_k) = \infty$,
then for every nonnegative integer $n$ we have $a_{k,n} = 0$
for all sufficiently large $k$.
The coefficient of $q^n$ in the partial sum $s_m(q)$ is
\beq           \label{qps:coeff}
\sum_{k=0}^{m} a_{k,n}.
\eeq
Since this sum is constant for all sufficiently
large $m$, it follows that the sequence $\{s_m(q)\}_{m=0}^{\infty}$ converges.

Conversely, if $\lim_{k\rightarrow\infty} v(f_k) \neq \infty$,
then there exists a nonnegative integer $n$ such that $a_{k,n} \neq 0$
for infinitely many integers $k$.
It follows that the sum~(\ref{qps:coeff}) is {\em not} constant for sufficiently large $m$,
and so the sequence $\{s_m(q)\}_{m=0}^{\infty}$ does not converge.
This completes the proof.
\eop

Theorem~\ref{qps:theorem:infiniteseries} allows us to substitute one formal power series
into another.  Let $f(q)$ and $h(q)$ be formal power series with $v(h) \geq 1.$
If 
\[
f(q) = \sum_{n=0}^{\infty} a_nq^n,
\]
then we define the composite function
\beq              \label{qps:comp}
(f\circ h) (q) = f(h(q)) = \sum_{n=0}^{\infty} a_nh(q)^n.
\eeq
Since
\[
v(a_nh(q)^n) = v(a_n)+nv(h) \geq n,
\]
it follows that the infinite series~(\ref{qps:comp}) converges.
Moreover,
\beq           \label{qps:vcompos}
v(f\circ h) = v(f)v(h).
\eeq

\bt               \label{qps:theorem:mkconvzero}
Let $a_0 \in R.$
If $\{f_k(q)\}_{k=1}^{\infty}$ is a sequence of formal power series
such that $f_k(q)$ has constant term $a_0$ for all $k \geq 1$, and 
if $\{m_k\}_{k=1}^{\infty}$ is a sequence of nonnegative integers such that
\[
\lim_{k\rightarrow\infty} m_k = \infty,
\]
then
\beq        \label{qps:mkconv}
\lim_{k\rightarrow\infty} f_k(q^{m_k}) = a_0.
\eeq
In particular, for any formal power series $f(q)$, 
\beq        \label{qps:miconv}
\lim_{k\rightarrow\infty} f(q^{m_k}) = f(0).
\eeq
\et

\pf
We have $v(f_k-a_0) \geq 1$ for all $k \geq 1.$
The convergence of~(\ref{qps:mkconv}) follows from~(\ref{qps:vcompos}) 
and the observation that
\[
v\left( f_k(q^{m^k})-a_0\right)
= v\left( (f_k-a_0)(q^{m_k}) \right)
= v\left( f_k-a_0\right)m_k
\geq m_k,
\]
hence
\[
\lim_{k\rightarrow\infty} v\left( f_k(q^{m^k})-a_0\right) = \infty.
\]
The convergence of~(\ref{qps:miconv}) follows from Theorem~\ref{qps:theorem:vconverge}.
\eop

\bt              \label{qps:theorem:kl}
Let $(k_i,\ell_i)_{i=0}^{\infty}$ be a sequence of ordered pairs of nonnegative integers
such that every ordered pair $(k,\ell)$ occurs exactly once in the sequence.
If the infinite series
\[
f(q) = \sum_{k=0}^{\infty} f_k(q)
\]
and
\[
g(q) = \sum_{\ell =0}^{\infty} g_{\ell}(q)
\]
converge in $R[[q]]$, then the infinite series
\beq      \label{qps:seriesi}
\sum_{i=0}^{\infty} f_{k_i}(q)g_{\ell_i}(q)
\eeq
converges, and
\[
f(q)g(q) = \sum_{i=0}^{\infty} f_{k_i}(q)g_{\ell_i}(q).
\]
In particular,
\beq     \label{qps:multseries}
f(q)g(q) = \sum_{n=0}^{\infty} \sum_{k+\ell=n} f_{k}(q)g_{\ell}(q)
\eeq
\et

\pf
By Theorem~\ref{qps:theorem:infiniteseries},
\[
\lim_{k\rightarrow\infty} v(f_k)
= \lim_{\ell\rightarrow\infty} v(g_{\ell}) = \infty,
\]
and so, for every $N$, there exists an integer $j_0 = j_0(N)$
such that $v(f_k) \geq N $ and $v(g_{\ell}) \geq N $ 
for all $k > j_0$ and $\ell > j_0$.
There are only finitely many ordered pairs $(k_i,\ell_i)$ with
both $k_i \leq j_0$ and $\ell_i \leq j_0$,
and so there is an integer $i_0 = i_0(N)$
such that, for all $i > i_0$, we have
$k_i > j_0$ or $\ell_i > j_0$, hence 
\[
v(f_{k_i}g_{\ell_i}) = v(f_{k_i})+v(g_{\ell_i}) \geq N
\]
and
\[
\lim_{i\rightarrow\infty} v(f_{k_i}g_{\ell_i}) = \infty.
\]
By Theorem~\ref{qps:theorem:infiniteseries},
the series~(\ref{qps:seriesi}) converges.
Moreover, for $m \geq i_0$,
\begin{align*}
s_m(q) & = \sum_{i=0}^{m} f_{k_i}(q)g_{\ell_i}(q) \\
& \equiv \sum_{i=0}^{i_0} f_{k_i}(q)g_{\ell_i}(q)  \pmod{q^N} \\
& \equiv \sum_{k=0}^{j_0} f_k(q)\sum_{\ell=0}^{j_0}g_{\ell}(q) \pmod{q^N} \\
& \equiv f(q)g(q) \pmod{q^N},
\end{align*}
and so
\[
v(fg-s_m) \geq N
\]
and
\[
\lim_{m\rightarrow\infty} v(fg-s_m) = \infty.
\]
By Theorem~\ref{qps:theorem:vconverge},
\[
\sum_{i=0}^{\infty} f_{k_i}(q)g_{\ell_i}(q) 
= \lim_{m\rightarrow\infty} s_m(q) = f(q)g(q).
\]
Identity~(\ref{qps:multseries}) is a special case of this result.
\eop

We define convergence of the {\em infinite product} \index{infinite product}
\[
\prod_{k=1}^{\infty} f_k(q)
\]
of formal power series as follows.
The {\em $m$--th partial product}\index{partial product}
\[
p_m(q) = \prod_{k=1}^{m} f_k(q)
\]
is simply a finite product of formal power series.
If the sequence $\{p_m(q)\}_{m=1}^{\infty}$ converges to $f(q) \in R[[q]]$,
then we write that the infinite product converges to $f(q)$, that is,
\[
\prod_{k=1}^{\infty} f_k(q) = \lim_{m\rightarrow\infty} p_m(q) = f(q).
\]
For example, from the unique representation of an integer
as the sum of distinct powers of 2, we have
\[
\prod_{k=1}^{m} (1+q^{2^{k-1}}) = (1+q)(1+q^2)(1+q^4)\cdots (1+q^{2^{m-1}})
= \sum_{n=0}^{2^m-1}q^n,
\]
and so
\[
\prod_{m=1}^{\infty} (1+q^{2^{m-1}}) = \sum_{n=0}^{\infty} q^n = (1-q)^{-1}.
\]

\bt                    \label{qps:theorem:product}
Let $\{f_k(q)\}_{k=1}^{\infty}$ be a sequence of formal power series
such that $f_k(q)$ has constant term $f_k(q)$ 1 for all $k \geq 1.$
The infinite product $\prod_{k=1}^{\infty} f_k(q)$ converges
if and only if $\lim_{k\rightarrow\infty} v(f_k(q)-1) = \infty$
if and only if the infinite series $\sum_{k=1}^{\infty} (f_k(q)-1)$ 
converges.
\et

\pf
This follows immediately from the definition of convergence.
\eop

\subsection{Change of variables in formal power series}

The next result allows us to substitute variables in identities
of formal power series.

\bt                           \label{qps:theorem:convergence1}
Let $\{g_k(q)\}_{k=1}^{\infty}$ be a sequence of formal power series
in $R[[q]]$ that converges to $g(q)$
and let $\{h_{\ell}(q)\}_{{\ell}=1}^{\infty}$ be a sequence
of formal power series that converges to $h(q)$.
If
\[
v(h_{\ell}) \geq 1 \qquad\text{for all $\ell \geq 1$},
\]
then
\[
\lim_{k,\ell\rightarrow\infty} g_k(h_{\ell}(q)) = g(h(q)).
\]
\et

\pf
For every nonnegative integer $N$ and formal power series
\[
f(q) = \sum_{n=0}^{\infty} a_nq^n,
\]
let
\[
f^{(N)}(q) = \sum_{n=0}^{N-1} a_nq^n.
\]
Then
\[
f(q) \equiv f^{(N)}(q) \pmod{q^{N}}.
\]
Since $v(h_{\ell}) \geq 1$ for all $\ell$, it follows that $v(h) \geq 1$,
and so the formal power series  $g(h(q))$ and $g_k(h_{\ell}(q))$
are well-defined for all nonnegative integers $k$ and $\ell.$
Since $\lim_{k\rightarrow\infty} g_k(q) = g(q)$, there exists an integer
$k_0(N)$ such that
\[
g_k^{(N)}(q) = g^{(N)}(q) \qquad\text{for all $k \geq k_0(N)$.}
\]
Similarly, since $\lim_{\ell\rightarrow\infty} h_{\ell}(q) = h(q)$,
there exists an integer ${\ell}_0(N)$ such that
\[
h_{\ell}^{(N)}(q) = h^{(N)}(q) \qquad\text{for all $\ell \geq {\ell}_0(N)$.}
\]
Therefore, for $k \geq k_0(N)$ and $\ell \geq \ell_0(N)$ we have
\[
g_k\left(h_{\ell}(q)\right) \equiv g_k^{(N)}\left(h_{\ell}^{(N)}(q)\right)
= g^{(N)}\left(h^{(N)}(q)\right) \equiv g(h(q)) \pmod{q^N}.
\]
This completes the proof.
\eop

\bt                           \label{qps:theorem:convergence2}
Let $\{g_k(q)\}_{k=1}^{\infty}$ be a sequence of formal power series
in $R[[q]]$ that converges to $g(q) \in R[[q]]$,
and let $h(q)$ be a formal power series with $v(h) \geq 1$.
Then
\[
\lim_{k\rightarrow\infty} g_k(h(q)) = g(h(q)).
\]
\et

\pf
This follows from Theorem~\ref{qps:theorem:convergence1}
with $h_{\ell}(q) = h(q)$ for all $\ell \geq 1$.
\eop

Theorem~\ref{qps:theorem:convergence2} is often applied in the following form.

\bt         \label{qps:theorem:sumproduct}
Let $\{f_k(q)\}_{k=0}^{\infty}$ and  $\{g_k(q)\}_{k=1}^{\infty}$
be sequences of formal power series such that
the infinite series $\sum_{k=1}^{\infty}f_k(q)$ converges
and the infinite product $\prod_{k=1}^{\infty}g_k(q)$ converges.
Let $h(q)$ be a formal power series with $v(h) \geq 1$.
If
\[
\sum_{k=0}^{\infty}f_k(q) = \prod_{k=1}^{\infty}g_k(q),
\]
then
\[
\sum_{k=0}^{\infty}f_k(h(q)) = \prod_{k=1}^{\infty}g_k(h(q)).
\]
\et

\pf
Let
\[
\sum_{k=0}^{\infty}f_k(q)= \sum_{k=0}^{m-1} f_k(q) = F(q).
\]
Theorem~\ref{qps:theorem:convergence2} implies that
\[
\sum_{k=0}^{\infty}f_k(h(q)) =
\lim_{m\rightarrow\infty} \sum_{k=0}^{m-1} f_k((h(q)) = F(h(q)).
\]
Similary,
\[
\prod_{k=1}^{\infty}g_k(q) =
\lim_{m\rightarrow\infty} \prod_{k=1}^{m}g_k(q) = F(q)
\]
implies that
\[
\prod_{k=1}^{\infty}g_k(h(q)) =
\lim_{m\rightarrow\infty} \prod_{k=1}^{m}g_k(h(q)) = F(h(q)).
\]
It follows that
\[
\sum_{k=1}^{\infty}f_k(h(q)) = \prod_{k=1}^{\infty}g_k(h(q)).
\]
This completes the proof.
\eop

\bt                           \label{qps:theorem:convergence3}
Let $\{g_k(q)\}_{k=1}^{\infty}$ be a sequence of formal power series
in $R[[q]]$, and let $h(q)$ be a formal power series with $v(h) \geq 1$.
If
\[
\lim_{k\rightarrow\infty} g_k(h(q)) = g(h(q)),
\]
then
\[
\lim_{k\rightarrow\infty} g_k(q) = g(q).
\]
\et

\pf
Let
\[
g_k(q) = \sum_{n=0}^{\infty} b_{k,n}q^n
\]
and
\[
g(q) = \sum_{n=0}^{\infty} b_{n}q^n.
\]
We must show that for every nonnegative integer $n$
we have $b_{k,n} = b_n$ for all sufficiently large $k$.
The proof will be by induction on $n$.
Let
\[
h(q) = \sum_{m=M}^{\infty} c_{m}q^m,
\qquad\text{where $v(h) = M \geq 1$ and $c_M \neq 0$.}
\]
Then
\[
h(q) \equiv c_Mq^M \pmod{q^{M+1}}.
\]
The congruences
\[
g_k(h(q)) \equiv b_{k,0} \pmod{q}
\]
and
\[
g(h(q)) \equiv b_{0} \pmod{q}
\]
imply that $b_{k,0} = b_0$ for all sufficiently large $k$.
Let $N \geq 1$.  Suppose there exists an integer $k_0(N)$ such that,
for $k \geq k_0(N)$ and $0 \leq n \leq N-1$, we have
$b_{k,n} = b_n$.
Then
\[
\sum_{n=N}^{\infty} b_{k,n}(h(q))^n \equiv b_{k,N}\left( c_M q^{M} \right)^{N}
\equiv b_{k,N}c_M^N q^{MN} \pmod{q^{MN+1}}
\]
and
\[
\sum_{n=N}^{\infty} b_{n}(h(q))^n \equiv b_{N}\left( c_M q^{M} \right)^{N}
\equiv b_{N}c_M^N q^{MN} \pmod{q^{MN+1}}.
\]
Since
\[
\lim_{k\rightarrow\infty} \sum_{n=N}^{\infty} b_{k,n}(h(q))^n
= \sum_{n=N}^{\infty} b_{n}(h(q))^n,
\]
it follows that for suffficiently large $k$ we have
\[
b_{k,N} c_M^{MN}= b_Nc_M^{MN},
\]
and so
\[
b_{k,N} = b_N.
\]
This completes the proof.
\eop

\bt                      \label{qps:theorem:substitution}
Let $R$ and $S$ be rings, and let $\varphi:R \rightarrow S[[q]]$
be a ring homomorphism.  Let $v$ be the valuation in the formal power series ring
$S[[q]]$.  Let $h(q) \in S[[q]]$ with $v(h) \geq 1$.
The map $\Phi: R[[q]] \rightarrow S[[q]]$ defined by
\beq        \label{qps:compo}
\Phi\left( \sum_{n=0}^{\infty} a_n q^n \right)
= \sum_{n=0}^{\infty} \varphi(a_n) h(q)^n
\eeq
is a ring homomorphism.
\et

\pf
Since
\[
v\left( \varphi(a_n) h(q)^n \right) =
v\left( \varphi(a_n) \right) + nv\left( h \right) \geq nv(h) \geq n,
\]
it follows from Theorem~\ref{qps:theorem:infiniteseries} that the infinite series
on the right side of~(\ref{qps:compo}) converges in the
formal power series ring $S[[q]]$.
Let $f(q) = \sum_{n=0}^{\infty}a_nq^n $
and $g(q) = \sum_{n=0}^{\infty} b_nq^n$.
Then
\begin{align*}
\Phi(f(q)) + \Phi(g(q))
& =  \sum_{n=0}^{\infty} \varphi(a_n) h(q)^n
+ \sum_{n=0}^{\infty} \varphi(b_n) h(q)^n \\
& =  \sum_{n=0}^{\infty} \left( \varphi(a_n) + \varphi(b_n) \right) h(q)^n \\
& =  \sum_{n=0}^{\infty} \varphi \left( a_n + b_n \right) h(q)^n \\
& =  \Phi \left( \sum_{n=0}^{\infty}  \left( a_n + b_n \right) q^n \right) \\
& =  \Phi \left( f(q) + g(q)\right).
\end{align*}
Similarly,
\begin{align*}
\Phi(f(q)) \Phi(g(q))
& =  \left( \sum_{n=0}^{\infty} \varphi(a_n) h(q)^n \right)
\left( \sum_{n=0}^{\infty} \varphi(b_n) h(q)^n \right) \\
& =  \sum_{n=0}^{\infty}
\left(\sum_{k+\ell = n}\varphi(a_k) \varphi(b_{\ell}) \right) h(q)^n \\
& =  \sum_{n=0}^{\infty}
\varphi\left(\sum_{k+\ell = n} a_k b_{\ell} \right) h(q)^n \\
& =  \Phi\left( \sum_{n=0}^{\infty}
\left(\sum_{k+\ell = n} a_k b_{\ell} \right) q^n \right) \\
& =  \Phi\left( \left( \sum_{n=0}^{\infty} a_n q^n \right)
\left( \sum_{n=0}^{\infty} b_n q^n \right)\right) \\
& =  \Phi \left( f(q)g(q)\right).
\end{align*}
This completes the proof.
\eop

\section{Limits of solutions of the functional equation for quantum multiplication}

The main result (Theorem~\ref{qps:theorem:infinitePlimit}) in this section 
states that if \pol\ is a solution of \tfe\
with $\supp(\mathcal{F}) = S(P)$ for a nonempty set $P$ of prime numbers,
and if $f_p(0) = 1$ for all $p \in P,$
then there exists a formal power series $F(q)$ such that
\[
F(q) = \lim_{\substack{n\rightarrow\infty \\ n \in S(P)}}f_n(q).
\]
For example, the sequence $\mathcal{F} = \{[n]_q\}_{n=1}^{\infty}$ has support 
$\supp(\mathcal{F}) = \N = S(P),$ where $P$ is the set of all primes,
and $\lim_{n\rightarrow\infty} [n]_q = (1-q)^{-1}.$

\bt                     \label{qps:theorem:ai}
Let \pol\ be a sequence of polynomials that satisfies 
the functional equation
\[
f_{mn}(q) = f_m(q)f_n(q^m).
\]
Let $\{a_i\}_{i=1}^{\infty}$ be a sequence of positive integers
such that $a_i \geq 2$ for infinitely many $i.$  
Let $n_0 = 1$ and
\[
n_k = a_1a_2\cdots a_k
\]
for $k \geq 1.$
If $f_{a_i}(0) = 1$ for all $i \geq 1,$ then the infinite product
\[
\prod_{i=1}^{\infty} f_{a_i}(q^{n_{i-1}})
\]
converges, and
\beq                \label{qps:qprod}
\lim_{k\rightarrow\infty} f_{n_k}(q) = \prod_{i=1}^{\infty}f_{a_i}(q^{n_{i-1}}).
\eeq
In particular, if $a_i = a \geq 2$ for all $i$, then there exists a  
formal power series $F_{a^{\infty}}(q)$ such that
\beq             \label{qps:ainfinity}
F_{a^{\infty}}(q) = \lim_{k\rightarrow\infty} f_{a^k}(q).
\eeq
\et

\pf
We shall show by induction that
\beq                       \label{qps:feprod}
f_{n_k}(q) = \prod_{i=1}^k f_{a_i}\left(q^{n_{i-1}}\right)
\eeq
for all $k\geq 1.$
The functional equation implies that
\begin{align*}
f_{n_1}(q) & = f_{a_1}(q) = f_{a_1}(q^{n_0}) \\
f_{n_2}(q) & = f_{a_1a_2}(q) = f_{a_1}(q)f_{a_2}(q^{a_1}) 
= f_{a_1}(q^{n_0})f_{a_2}(q^{n_1}) \\
f_{n_3}(q) & = f_{a_1a_2a_3}(q) = f_{a_1a_2}(q)f_{a_3}(q^{a_1a_2}) \\
& = f_{a_1}(q^{n_0})f_{a_2}(q^{n_1})f_{a_3}(q^{n_3}),
\end{align*}
and so formula~(\ref{qps:feprod}) holds for $k = 1,2,3$.
If~(\ref{qps:feprod}) holds for $k$, then
\begin{align*}
f_{n_{k+1}}(q) & = f_{n_ka_{k+1}}(q) \\
& = f_{n_k}(q)f_{a_{k+1}}(q^{n_k}) \\
& = \left(\prod_{i=1}^k f_{a_i}\left(q^{n_{i-1}}\right)\right)f_{a_{k+1}}(q^{n_k}) \\
& = \prod_{i=1}^{k+1} f_{a_i}\left(q^{n_{i-1}}\right).
\end{align*}
This completes the induction.

Since $f_{a_i}(0) = 1$ for all $i \geq 1,$ 
it follows that $f_{n_k}(0) = 1$ 
and $v(f_{n_k}(q)-1)\geq 1$ for all $k \geq 1.$
By Theorem~\ref{qps:theorem:mkconvzero},
\[
v\left(f_{n_k}\left(q^{n_{k-1}}\right)-1\right) \geq n_{k-1}.
\]
Since $\lim_{k\rightarrow\infty} n_k = \infty,$ 
Theorem~\ref{qps:theorem:product} implies the convergence 
of the infinite product~(\ref{qps:qprod}).
This immediately implies~(\ref{qps:ainfinity}).
\eop

\bt             \label{qps:theorem:akbk}
Let \pol\ be a sequence of polynomials that satisfies 
the functional equation
\[
f_{mn}(q) = f_m(q)f_n(q^m).
\]
If $a \geq 2$, $b \geq 2$, and $f_a(0) = f_b(0) = 1,$ then
\[
F_{a^{\infty}}(q) = \lim_{k\rightarrow\infty} f_{a^k}(q)
= \lim_{k\rightarrow\infty} f_{b^k}(q)= F_{b^{\infty}}(q).
\]
\et

\pf
The functional equation implies that $f_{a^k}(0) = f_{b^k}(0) = f_{(ab)^k}(0) = 1$
for all $k \geq 1.$  By Theorem~\ref{qps:theorem:mkconvzero},
\[
F_{b^{\infty}}(0) = \lim_{k\rightarrow\infty} f_{b^k}(q^{a^k}) = 1.
\] 
It follows from Theorem~\ref{qps:theorem:ai} that
\begin{align*}
F_{(ab)^{\infty}}(q) 
& = \lim_{k\rightarrow\infty} f_{(ab)^k}(q) \\
& = \lim_{k\rightarrow\infty} f_{a^k}(q)f_{b^k}(q^{a^k}) \\
& = \lim_{k\rightarrow\infty} f_{a^k}(q) \lim_{k\rightarrow\infty} f_{b^k}(q^{a^k}) \\
& = F_{a^{\infty}}(q). 
\end{align*}
Similarly, $F_{(ab)^{\infty}}(q) = F_{b^{\infty}}(q)$ 
and so $F_{a^{\infty}}(q) = F_{b^{\infty}}(q)$.
This completes the proof.
\eop

\bt                     \label{qps:theorem:finitePlimit}
Let \pol\ be a sequence of polynomials that satisfies 
the functional equation
\[
f_{mn}(q) = f_m(q)f_n(q^m).
\]
If $\supp(\mathcal{F}) = S(P)$ for some nonempty finite set $P$ of prime numbers,
and if $f_p(0)=1$ for all $p \in P,$
then there exists a formal power series $F(q)$ such that
\beq          \label{qps:finitePlimit}
F(q) = \lim_{\substack{n\rightarrow\infty \\ n \in S(P)}}f_n(q).
\eeq
\et

\pf
Since $f_p(0)=1$ for all $p \in P,$ it follows that $f_m(0)=1$,
$v(f_m) = 0$, and $v(f_m-1)\geq 1$ for all $m \in S(P).$

By Theorem~\ref{qps:theorem:akbk}, there exists a formal power series $F(q)$
such that
\[
F(q) = F_{p^{\infty}}(q)
= \lim_{k\rightarrow\infty} f_{p^k}(q) \qquad\text{for all $p \in P$.}
\]
Since the set $P$ is finite, for every integer $N$ there exists an integer 
\[
k_0=k_0(N) \geq N
\]
such that 
\[
v(F-f_{p^k}) \geq N \qquad\text{for all $p \in P$ and $k \geq k_0.$}
\]
Let $n_0 = n_0(N) = \prod_{p\in P} p^{k_0-1}$.
If $n > n_0$ and $n \in S(P)$, then $n = p^{k_0} m$
for some $p \in P$ and $m\in S(P)$.
Moreover,
\begin{align*}
F(q) - f_n(q)
& = F(q) - f_{p^{k_0}}(q) f_m(q^{p^{k_0}}) \\
& = F(q) - f_{p^{k_0}}(q) - f_{p^{k_0}}(q)\left( f_m(q^{p^{k_0}}) - 1\right).
\end{align*}
We have 
\[
v(F-f_{p^{k_0}}) \geq N
\]
and
\begin{align*}
v\left( f_{p^{k_0}}(q) \left( f_m(q^{p^{k_0}}) -1\right)  \right)
& = v\left( f_{p^{k_0}}\right) + v\left( f_m(q^{p^{k_0}}) -1 \right) \\
& = v\left( f_m(q^{p^{k_0}}) -1 \right) \\
& = v(f_m-1) p^{k_0} \\
& \geq p^{k_0}\\
& \geq N.
\end{align*}
It follows that
\[
v(F-f_n) \geq 
\min\left( v(F-f_{p^{k_0}}),v\left( f_{p^{k_0}}(q)\left( f_m(q^{p^{k_0}})-1\right)\right)\right)
\geq N
\]
for all $n \geq n_0.$ 
This completes the proof of~(\ref{qps:finitePlimit}).
\eop

\bt                     \label{qps:theorem:infiniteP}
Let \pol\ be a sequence of polynomials that satisfies 
the functional equation
\[
f_{mn}(q) = f_m(q)f_n(q^m).
\]
If $\supp(\mathcal{F}) = S(P)$, where $P$ is an infinite set of prime numbers,
and if $f_p(0)=1$ for all $p \in P,$
then there exists a formal power series $F_P(q)$ such that
\beq          \label{qps:infiniteP}
F_P(q) = \lim_{\substack{p\rightarrow\infty \\ p \in P}}f_p(q).
\eeq
\et

\pf
If $p_1$ and $p_2$ are prime numbers in $P$ such that $p_1 \geq N$ and $p_2 \geq N,$
then
\begin{align*}
f_{p_1}(q) - f_{p_2}(q)
& = f_{p_1}(q) - f_{p_1p_2}(q) + f_{p_1p_2}(q) - f_{p_2}(q)  \\
& = f_{p_1}(q) - f_{p_1}(q)f_{p_2}(q^{p_1}) + f_{p_2}(q)f_{p_1}(q^{p_2}) - f_{p_2}(q)  \\
& = f_{p_1}(q)\left( 1 - f_{p_2}(q^{p_1})\right) + 
    f_{p_2}(q)\left( f_{p_1}(q^{p_2}) - 1 \right).
\end{align*}
Since 
\[
v(f_{p_1}) = v(f_{p_2}) = 0
\]
and, by Theorem~\ref{qps:theorem:val},
\[
v\left( 1 - f_{p_2}(q^{p_1}) \right) \geq p_1 \geq N,
\]
and
\[
v\left( f_{p_1}(q^{p_2}) - 1 \right) \geq p_2 \geq N,
\]
it follows that
\[
v(f_{p_1} - f_{p_2}) \geq N.
\]
Therefore, 
\[
\lim_{p_1,p_2\rightarrow\infty} v(f_{p_1} -f_{p_2}) = \infty,
\]
and Theorem~\ref{qps:theorem:vconverge} implies the existence 
of the limit~(\ref{qps:infiniteP}).
This completes the proof.
\eop

\bt                     \label{qps:theorem:infinitePlimit}
Let \pol\ be a sequence of polynomials that satisfies 
the functional equation
\[
f_{mn}(q) = f_m(q)f_n(q^m).
\]
Let $\supp(\mathcal{F}) = S(P),$ where $P$ is a nonempty set of prime numbers.
If $f_p(0)=1$ for all $p \in P,$
then there exists a formal power series $F(q)$ such that
\beq          \label{qps:infinitePlimit}
F(q) = \lim_{\substack{n\rightarrow\infty \\ n \in S(P)}}f_n(q).
\eeq
\et

\pf
If $P$ is finite, then the result is Theorem~\ref{qps:theorem:finitePlimit}.

Suppose that $P$ is an infinite set of prime numbers.
Let $N$ be a positive integer.
By Theorem~\ref{qps:theorem:infiniteP}, there is a formal power series $F_P(q)$
such that 
\[
F_P(q) = \lim_{\substack{p\rightarrow\infty \\ p \in P}}f_p(q),
\]
and so there exists a prime number
\[
p_0 = p_0(N) \geq N
\]
such that
\[
v(F_P - f_p) \geq N \qquad\text{for all $p \in P$ with $p \geq p_0$.}
\]
If $n = pm \in S(P),$ where $p \in P$ and $p \geq p_0$, then
\begin{align*}
F_P(q) - f_n(q) 
& = F_P(q) - f_p(q)f_m(q^p) \\
& = \left( F_P(q) - f_p(q)\right) + f_p(q)\left( 1 - f_m(q^p)\right).
\end{align*}
Since $v\left( 1 - f_m(q^p)\right)\geq p > p_0 \geq N,$
it follows that $v(F_P - f_n) \geq N.$

By Theorem~\ref{qps:theorem:akbk}, there is a formal power series $F(q)$
such that
\[
F(q) = \lim_{k\rightarrow\infty} f_{p^k}(q)
\]
for all $p \in P.$  There exists an integer
\[
k_0 = k_0(N) \geq N
\]
such that
\[
v\left( F-f_{p^k}\right) \geq N
\qquad\text{for all $k \geq k_0$ and $p \in P$ with $p < p_0$.}
\]
If $n = p^km \in S(P),$ where $k \geq k_0$ and $p \in P$ with $p < p_0$, then
\begin{align*}
F(q) - f_n(q) 
& = F(q) - f_{p^k}(q)f_m(q^{p^k}) \\
& = \left( F(q) - f_{p^k}(q)\right) + f_{p^k}(q)\left( 1 - f_m(q^{p^k})\right).
\end{align*}
Since $v\left( 1 - f_m(q^{p^k})\right)\geq p^k \geq k \geq k_0 \geq N,$
it follows that $v(F - f_n) \geq N.$

Let $p_1$ and $p_2$ be prime numbers in $P$ such that $p_1 < p_0 \leq p_2,$
and let $k \geq k_0.$  Let $n = p_1^kp_2.$  Then
\[
v(F_P-f_n) \geq N,
\]
\[
v(F-f_n) \geq N,
\]
and so
\[
v(F-F_P) = v\left( (F-f_n) - (f_n-F_P) \right) \geq N.
\]
This holds for all $N$, and so 
\[
F(q) = F_P(q).
\]

Let $n_0 = \prod_{\substack{p<p_0 \\p \in P} } p^{k_0-1} = n_0(N)$ and 
let $n \in S(P)$ with $n \geq n_0.$
If $p$ divides $n$ for some prime $p \in P$, $p \geq p_0,$
then $v(F-f_n) = v(F_P-f_n) \geq N.$
If every prime factor of $n$ is less than $p_0$,
then $n$ is an integer of the form 
$n = \prod_{\substack{p<p_0 \\p \in P} } p^{a_p}$,
where $a_p \in \N_0$ for all $p \in P$. 
Since $n \geq n_0,$ it follows that $a_p \geq k_0$ for some $p$,
and so $v(F-f_n) \geq N.$  Therefore, 
\[
\lim_{\substack{n\rightarrow\infty \\ n \in S(P)}} v(F-f_n) = \infty.
\]
This completes the proof.
\eop

Let \pol\ be a sequence of polynomials satisfying \tfe.
We assume that $\supp(\mathcal{F}) = S(P),$ 
where $P$ is a nonempty set of prime numbers,
and $f_p(0)=1$ for all $p \in P.$
By Theorem~\ref{qps:theorem:infinitePlimit},
there exists a formal power series $F(q)$ such that
\[
F(q) = \lim_{\substack{n\rightarrow\infty \\ n \in S(P)}}f_n(q).
\]

\bp
Which formal power series are limits of solutions of \tfe?
\ep

\bp
To what extent is it possible to reconstruct the sequence $\mathcal{F}$
from the formal power series $F(q)$?
If $F(q)$ is the limit of a solution of \tfe, determine all solutions 
$\mathcal{F}$ of~(\ref{qps:fe}) whose limit is $F(q).$
\ep

Let \polg\ be a solution of \tfe\ with support $S(P_G)$ and with limit
\[
F(q) = \lim_{\substack{n\rightarrow\infty \\ n \in S(P_G)}}g_n(q).
\]
We observe that if $P$ is a nonempty subset of $P_G$, then we can construct 
another solution \pol\ of~(\ref{qps:fe}) by {\em restriction} 
of $\mathcal{G}$ to the subsemigroup 
$S(P)$ of $S(P_G)$ as follows:  
\[
f_n(q) = \left\{
\ba{ll}
g_n(q) & \text{if $n \in S(P)$,} \\
0      & \text{if $n \not\in S(P).$} 
\ea
\right.
\]
Since
\[
F(q) = \lim_{\substack{n\rightarrow\infty \\ n \in S(P)}}f_n(q),
\]
it follows that $F(q)$ does not uniquely determine $\mathcal{F}$.

The suggests the following definition.  A solution \pol\ of \tfe\
will be called {\em maximal} if $\mathcal{F}$ cannot be constructed 
from another solution $\mathcal{G}$ of the functional equation 
by restriction to a subsemigroup of $\supp(\mathcal{G}).$

\bt
Every solution of \tfe\ is contained in at least one maximal solution.
\et

\pf
Let \pol\ be a solution of~(\ref{qps:fe}) with $\supp(\mathcal{F}) = S(P).$
Let $X$ be the set of all pairs $(\mathcal{G},P_G)$ such that
\benum
\item[(i)]
$P_G$ is a set of prime numbers and $P \subseteq P_G$,
\item[(ii)]
$\mathcal{G}$ is a solution of \tfe\ with support $S(P_G)$,
\item[(iii)]
The restriction of $\mathcal{G}$ to $S(P)$ is $\mathcal{F}$.
\eenum
The set $X$ is partially ordered as follows:
$(\mathcal{G},P_G) \preceq (\mathcal{H},P_H)$ if $P_G \subseteq P_H$ 
and the restriction of $\mathcal{H}$ to $S(P_G)$ is $\mathcal{G}$.
Every chain in $X$ has an upper bound, and so, by Zorn's Lemma, 
the set $X$ contains a maximal element $(\mathcal{H},P_H)$.
Then $\mathcal{H}$ is a maximal solution of \tfe\ that restricts to $\mathcal{F}$.
\eop

\bp
Describe the maximal solutions of \tfe.  
For which sets $P$ of prime numbers does there exist a maximal solution 
$\mathcal{F}$ of~(\ref{qps:fe}) with $\supp(\mathcal{F}) = S(P)$?
\ep

\section{A functional equation for formal power series} \label{qps:section4}

We shall prove that the functional equation $f(q)F(q^m) = F(q)$ 
has a unique solution not only for polynomials $f(q)$ 
but also for formal power series $f(q)$ with constant term 1.

\bt
Let $m$ be a positive integer.
For every formal power series $f(q)$ with constant term $f(0) = 1,$
there exists a unique formal power series
\beq                           \label{qps:fpsoln}
F(q) = \prod_{i=1}^{\infty} f\left(q^{m^i}\right) = f(q) f(q^m) f\left(q^{m^2}\right)\cdots
\eeq
that satisfies the functional equation
\beq        \label{qps:fpeq2}
f(q)F(q^m) = F(q).
\eeq
\et

\pf
The infinite product~(\ref{qps:fpsoln}) converges in the ring of formal power series, since 
\[
v\left( f\left( q^{m^i} \right) -1 \right) =  v(f(q)-1)m^i
\]
and
\[
\lim_{i\rightarrow\infty} v\left(f\left(q^{m^i}\right)-1\right) = \infty.
\]
Then
\begin{align*}
f(q)F(q^m) 
& = f(q)\prod_{i=0}^{\infty} f\left(\left(q^m\right)^{m^i}\right)  \\
& = f(q)\prod_{i=0}^{\infty} f\left(q^{m^{i+1}}\right)  \\
& = f(q)\prod_{i=1}^{\infty} f\left(q^{m^{i}}\right)  \\
& = \prod_{i=0}^{\infty} f\left(q^{m^{i}}\right)  \\
& = F(q).
\end{align*}
This proves that the formal power series~(\ref{qps:fpsoln}) is a solution of~(\ref{qps:fpeq2}).

Conversely, let $F(q)$ be a formal power series that is a solution of~(\ref{qps:fpeq2}).
We shall prove that 
\beq         \label{qps:induc}
F(q) = \prod_{i=0}^{k-1} f\left(q^{m^{i}}\right)F\left(q^{m^k}\right)
\eeq
for every positive integer $k$.
From the functional equation~(\ref{qps:fpeq2}), we obtain
\[
F(q) = f(q)F(q^m) = f(q)\left( f(q^m)F((q^m)^m)\right)
 = f(q)f(q^m)F(q^{m^2}).
\]
This proves~(\ref{qps:induc}) for $k = 1$ and 2.
If~(\ref{qps:induc}) holds for some integer $k \geq 2$, then
\begin{align*}
F(q) 
& = f(q)\prod_{i=0}^{k-1} f\left(q^{m^{i+1}}\right)F\left(q^{m^{k+1}}\right)\\
& = \prod_{i=0}^{k} f\left(q^{m^{i}}\right)F\left(q^{m^{k+1}}\right),
\end{align*}
and so~(\ref{qps:induc}) holds for $k+1$.
This completes the induction.
Since $F(0)=1,$ it follows from Theorem~\ref{qps:theorem:mkconvzero} that
\begin{align*}
F(q) 
& = \lim_{k\rightarrow\infty}\prod_{i=0}^{k-1} f\left(q^{m^{i}}\right)
\lim_{k\rightarrow\infty}F\left(q^{m^k}\right) \\
& = \lim_{k\rightarrow\infty}\prod_{i=0}^{k-1} f\left(q^{m^{i}}\right) \\
& = \prod_{i=0}^{\infty} f\left(q^{m^{i}}\right).
\end{align*}
This completes the proof.
\eop

If $m$ is a positive integer and $f_m(q)$ a polynomial of degree less than $m$,
then we can explicitly construct the coefficients of the unique formal power series that
satisfies the equation $f_m(q)F(q^m) = m(q).$
We need the following notation.
Let $m$ be an integer, $m \geq 2.$  Every integer $k$ has
a unique $m$--adic representation
\[
k = \sum_{t=0}^{\infty} e_{t,k} m^{t},
\]
where
\[
 e_{t,k} \in \{0,1,2,\ldots, m-1\}
\]
for all $t,$ and $e_{t,k} = 0$ for all but finitely many $t.$
For $i = 1,2,\ldots,m-1,$ let 
\[
d_i(k) = \sum_{e_{t,k} = i}1.
\]
Then $d_i(k)$ is the arithmetic function that counts the number of 
digits equal to $i$ in the $m$--adic representation of $k$.
For example, if $m = 3,$ then $71 = 2\cdot 3^0 + 2\cdot 3^1 + 1\cdot 3^2+ 2\cdot 3^3,$
and so $d_1(71) = 1$ and $d_2(71)=3.$

\bt
Let $f_m(q) = 1 + \sum_{i=1}^{m-1}a_iq^i$ be a polynomial of degree $m-1$.
Let $F(q) = 1+ \sum_{k=1}^{\infty}b_k q^k$ be the formal power series with coefficients
\beq         \label{qps:digit}
b_k = \prod_{i=1}^{m-1} a_i^{d_i(k)}
\qquad\text{for all $k = 1,2,3,\ldots.$}
\eeq
Then
\beq            \label{qps:fepow}
f_m(q)F(q^m) = F(q).
\eeq
Conversely, if $F(q) = 1+ \sum_{k=1}^{\infty}b_k q^k$ 
is a formal power series that satisfies the functional equation~(\ref{qps:fepow}), 
then the coefficients of $F(q)$ are given by~(\ref{qps:digit}).
\et

\pf
Let $a_0 = b_0 = 1.$
From  the functional equation~(\ref{qps:fepow}) we have
\[
\left(\sum_{i=0}^{m-1}a_iq^i\right) \left( \sum_{j=0}^{\infty}b_j q^{mj} \right)
= \sum_{i=0}^{m-1} \sum_{j=0}^{\infty} a_ib_j q^{i+mj}
= \sum_{k=0}^{\infty} b_k q^{k},
\]
if and only if
\[
a_ib_j = b_{i+mj}
\]
for all $i = 0,1,\ldots,m-1$ and $j \in \N_0.$
We shall prove that the the unique solution of~(\ref{qps:psfe}) is~(\ref{qps:digit}).

For $j=0$, we have 
\[
b_k =a_k = \prod_{i=0}^{m-1} a_i^{d_i(k)}
\]
for $k=0,1,\ldots,m-1.$
This proves~(\ref{qps:digit}) for $0 \leq k < m.$

Suppose that~(\ref{qps:digit}) holds for all nonnegative integers $k < m^{r}.$
Let $k$ be an integer in the interval $m^r \leq k < m^{r+1}.$
We can write 
\[
k = i + mj, \qquad\text{where $0 \leq i \leq m-1$ 
and $m^{r-1} \leq j < m^r.$}
\]
If
\[
j = \sum_{t=0}^{\infty} e_{t,j} m^{t},
\]
then
\[
mj = \sum_{t=0}^{\infty} e_{t,j} m^{t +1}
= \sum_{t=1}^{\infty} e_{t,mj} m^{t},
\]
and
\[
e_{t,mj} = e_{t-1,j} 
\qquad\text{for  all $t \in \N.$}
\]
Then
\[
k = i + mj = i + \sum_{t=1}^{\infty} e_{t,mj} m^{t},
\]
and so
\[
d_{t}(k) = \left\{\ba{ll}
d_{t}(j) & \text{for $t \neq i$}  \\
d_{t}(j) + 1 & \text{for $t = i$.}  
\ea\right.
\] 
Therefore,
\begin{align*}
\prod_{t=0}^{m-1} a_t^{d_t(k)}
& =  a_i^{d_i(k)}\prod_{\substack{t=0 \\ t\neq i}}^{m-1} a_t^{d_t(k)} \\
& =  a_i^{d_i(j)+1}\prod_{\substack{t=0 \\t\neq i}}^{m-1} a_t^{d_t(j)} \\
& =  a_i\prod_{t=0}^{m-1} a_t^{d_t(j)} \\
& =  a_ib_j \\
& =  b_k.
\end{align*}
This completes the proof.
\eop

{\em Acknowledgement.}  
I wish to thank Yang Wang for a helpful discussion about self-similar sets.

\providecommand{\bysame}{\leavevmode\hbox to3em{\hrulefill}\thinspace}
\providecommand{\MR}{\relax\ifhmode\unskip\space\fi MR }
\providecommand{\MRhref}[2]{%
  \href{http://www.ams.org/mathscinet-getitem?mr=#1}{#2}
}
\providecommand{\href}[2]{#2}


\begin{thebibliography}{1}


\bibitem{nath03b}
M.~B. Nathanson, \emph{A functional equation arising from multiplication of quantum
  integers}, J. Number Theory, to appear.

\bibitem{nath03d}
\bysame, \emph{Additive number theory and the ring of quantum
  integers}, www.arXiv.org, math.NT/0204006.

\end{thebibliography}

\end{document}